\documentclass[12pt,rieqno,a4paper,twoside]{article}
\usepackage{latexsym}
\usepackage{amsfonts}
\usepackage{amsthm}
\usepackage{amsmath}
\usepackage[mathscr]{eucal}
\usepackage{epsfig}
\usepackage{a4wide}
\usepackage{ifthen}

\usepackage{graphicx}
\usepackage{tabularx}
\usepackage{fancyhdr}

\def\@copyright{}
\theoremstyle{plain}

\newtheorem{theorem}{Theorem}[section]

\newtheorem{example}{Example}[section]

\newtheorem{lemma}{Lemma}[section]

\newtheorem{corollary}{Corollary}[section]
\newtheorem{remark}{Remark}[section]

\newcommand{\N}{ \mathbb{N} }

\newcommand{\R}{ \mathbb{R} }

\newcommand{\trunc}[1]{ {\lfloor #1 \rfloor} }

\newcommand{\calE}{\mathcal{E}}

\newcommand{\calR}{\mathcal{R}}

\newcommand{\calV}{\mathcal{V}}

\newcommand{\eins}{{\mathbf 1}}
\newcommand{\matA}{{\mathbf A}}

\newcommand{\matH}{{\mathbf H}}

\newcommand{\matX}{{\mathbf X}}

\newcommand{\matW}{{\mathbf W}}

\newcommand{\veca}{{\mathbf a}}

\newcommand{\vecv}{{\mathbf v}}

\newcommand{\vecx}{{\mathbf x}}

\newcommand{\vecY}{{\mathbf Y}}

\newcommand{\bfbeta}{\boldsymbol{\beta}}

\newcommand{\bfeps}{\boldsymbol{\epsilon}}
\newcommand{\bhbeta}{\widehat{\boldsymbol{\beta}}}

\newcommand{\Cov}{{\mbox{Cov\,}}}

\newcommand{\diag}{{\mbox{diag\,}}}
\newcommand{\matid}{{\mathbf{I}}}

\newcommand{\id}{\operatorname{id}}

\numberwithin{equation}{section}

\begin{document}

\thispagestyle{empty}

\pagestyle{fancy}

\begin{center}
  \Large \bf
  Sequentially Updated Residuals and Detection of Stationary
  Errors in Polynomial Regression Models
\end{center}
\vskip 1cm
\begin{center}
  \textbf{Ansgar Steland}  \\
Institute of Statistics, RWTH Aachen University, Aachen, Germany
\end{center}

\hskip 2cm

\noindent {\bf Abstract:} The question whether a time series
behaves as a random walk or as a stationary process is an
important and delicate problem, particularly arising in financial
statistics, econometrics, and engineering. This paper studies the
problem to detect sequentially that the error terms in a
polynomial regression model no longer behave as a random walk but
as a stationary process. We provide the asymptotic distribution
theory for a monitoring procedure given by a control chart, i.e.,
a stopping time, which is related to a well known unit root test
statistic calculated from sequentially updated residuals. We
provide a functional central limit theorem for the corresponding
stochastic process which implies a central limit theorem for the
control chart. The finite sample properties are investigated by a
simulation study.



\vfill

\noindent
\quad Address correspondence to A. Steland, Institute of Statistics, RWTH Aachen University,
W\"ullnerstra\ss e 3, 52056 Aachen, Germany; Fax: +492418094573; E-mail:
steland@stochastik.\\rwth-aachen.de

\newpage
\section{INTRODUCTION}

\noindent Random walks have been proposed as reasonable models for
discretely observed data in many disciplines. In engineering, they
have been proposed to model production processes with degradation.
For instance, the additive damage model assumes that damage
cumulates yielding a random walk, and the system fails if the
cumulative damage reaches a threshold. We refer to Birnbaum and
Saunders (1969), Taguchi (1981, 1985), Taguchi et al. (1989), Adams
and Woodall (1989), Doksum and H\'oyland (1992), Vander Wiel (1996),
Durham and Padgett (1997), Park and Padgett (2006), and Srivastava
and Wu (1994, 2003).
In financial statistics, random walks appear as a model for the
(log) prices of an exchange-traded asset. That idea dates back to
Bachelier (1900), and nowadays there is an extensive literature on
the random walk hypothesis in the empirical finance literature,
mainly addressing the question whether the increments are
correlated. Random walks have also been proposed as a model for
important economic series as the gross domestic product.
Therefore, an important problem is to check sequentially whether a
time series is compatible with the random walk model or follows an
alternative (out-of-control) model under which the series is
stationary.

As is well known, a false answer to that question can lead to
completely wrong statistical conclusions, since even elementary
statistics change their convergence rates and limit distributions.
The implications for a rich class of nonparametric kernel control
charts covering, e.g., an approximation to the classic EWMA
control chart have been discussed in detail in Steland (2004).
Another popular approach to monitor both i.i.d. observations and
random walks resp. Brownian motions to detect changes in the mean
is based on the CUSUM procedure, which is known to be optimal in
the sense of Lorden's criterion. We refer to Beibel (1996),
Moustakides (1986, 2004, 2007), Ritov (1990), Siegmund (1985), Shiryaev (1996), and to
the monograph of Brodsky and Darkhovsky
(2000). Having this in mind, it is of particular interest to study
sequential monitoring (surveillance) procedures, which are
designed to detect departures from the random walk hypothesis as
soon as possible.

In this article we investigate a sequential monitoring procedure
which is related to a well known unit root test studied in detail by
Breitung (2002). To test the unit root null hypothesis against the
alternative of stationarity, he proposed to use a variance ratio
statistic comparing the dispersion of partial sums with the
dispersion of the observations. That test statistic is similar to
the statistic underlying the so-called variance ratio or KPSS test
proposed by Kwiatkowski et al. (1992) to test the inverse testing
problem of stationarity against the unit root alternative. Lee and
Schmidt (1996) have shown that the KPSS test is also consistent
against stationary long-memory alternatives, for a further detailed
study we refer to Giraitis et al. (2003). The KPSS test is known to
be powerful for many important data generating processes and robust
in terms of the type I error rate. For both testing problems (random
walk versus stationarity and vice versa) sequential monitoring
(surveillance) procedures based on control charts related to the
variance ratio statistic have been proposed in Steland (2007a). In
that paper the original time series $ Y_1, Y_2, \ldots $ is
monitored. Under mild conditions the asymptotic distributions of the
associated stopping times have been established under various
in-control and change-point models.

Motivated by promising results from a preliminary study (Steland,
2006), this article considers the more involved and delicate
problem to test sequentially whether or not the error terms in a
polynomial regression model form a random walk, thus allowing for
nonlinear time trends. Assume we observe sequentially a time
series $ \{ Y_t : t \in \N \} $ of real-valued observations
satisfying
\[
  Y_t = m_t + \epsilon_t, \qquad t \in \N,
\]
with $ E(\epsilon_t) = 0 $ for all $ t $. In many applications the
regression function $ m $ is smooth, which motivates to consider
polynomials of known degree. Thus, we assume
\begin{equation}
\label{ModelOfPaper}
  Y_t = \beta_0 + \beta_1 t + \dots + \beta_p t^p + \epsilon_t, \qquad t \in \N,
\end{equation}
where $ \bfbeta = (\beta_0,\dots,\beta_p)' \in \R^{p+1} $ are
unknown regression coefficients and $ p \in \N_0 $. Basically, the
aim is to detect a departure from the in-control model that the
error terms form a random walk in favor of a stationary process.
Note that the model covers the case that before the change a
Brownian motion with polynomial drift, $ \xi(t) = \mu(t) + \sigma
B(t) $, where $ B $ denotes standard Brownian motion, $ \sigma > 0
$ is a constant, and $ \mu(t) = \sum_{j=0}^p \beta_j t^j $,  is
discretely sampled at time instances $ t = 1, 2, \dots, q-1 $. In
this case $ \epsilon_t = \sigma B(t) \sim N( 0, \sigma^2 t ) $,
i.e., the variance is a linear function of time. After the change
we observe $ \xi(t) = \mu(t) + \sigma B(q)+\eta(t) $, $
t=q,q+1,\ldots, $ where $\eta(t)$ is a stationary process; e.g.
given by a continuous-time moving average,
\[
\eta(t) = \int^t_{-\infty} \varphi(t-s)\,dB(s),
\]
for some function $ \varphi $ with $ \int \varphi^2(t) \, dt <
\infty $. Our results allow for substantially more general error
sequences.

Since for many practical applications the most important
alternative model is a (polynomial) time trend with stationary
errors, we will apply a control chart (stopping time) providing a
signal, if there is evidence that the errors are no longer
compatible with the random walk hypothesis. We provide the
asymptotic distribution theory under the in-control model that the
error terms behave as a random walk but allow for an unknown
polynomial time trend. Further, we establish results under a
change-point model where the errors form a stationary process
after an unknown change-point. Since our results provide the
asymptotic distribution of the stopping time, one may design a
surveillance procedure according to various criteria.
Particularly, our results allow to design the procedure to
guarantee a specified asymptotic significance level (type I error
rate). If we get a signal, the classic polynomial regression model
with stationary errors can be regarded as statistically confirmed,
which is an attractive property for many applications.

We study the intuitive approach to calculate the least squares
residuals and to apply an appropriate monitoring procedure to
these residuals. In sequential analysis, recursive residuals are
often used, see the classic paper by Brown et al. (1975), and Sen (1982), mainly because they are fast to compute.
However, having in mind contemporary computing facilities, we
introduce {\em sequentially updated residuals}, where at each step
the full set of residuals is calculated. We consider a monitoring
procedure with a time horizon $ T $ where monitoring stops,
because in many real applications it is unrealistic to assume that
monitoring can be conducted forever. However, the modifications of
the results to allow for infinite monitoring are straightforward
and briefly discussed.

The rest of the paper is organized as follows. In
Section~\ref{SecModel} we specify and discuss the assumptions on
the error terms and introduce the proposed procedure and required
regularity conditions. The asymptotic results for the process of
sequentially updated residuals, the process associated to the
proposed control statistic, and for the resulting stopping time,
are discussed in detail in Section~\ref{SecIP} under the
in-control model that the regression errors behave as a random
walk. The results are constructive in the sense that explicit
representations of the asymptotic error process can be obtained in
terms of the moment functions, $  \int_0^s t^k B(t) \,dt $, $ s
\in [0,1] $, $k \in \N $, associated to a standard Brownian motion
$B$, which makes simulation from the limiting processes feasible.
Section~\ref{SecCP} gives asymptotic results under a change-point
model where the behavior changes after a certain fraction of the
data from a random walk behavior to a stationary process. We
report in Section~\ref{SecSim} about a simulation study which
examines some finite sample properties of the method. Proofs of
the main results of this paper are postponed to appendices.

\section{MODEL, ASSUMPTIONS, AND THE METHOD}
\label{SecModel}

\subsection{Model and Assumptions}

\noindent It remains to specify model (\ref{ModelOfPaper}) in
detail. We assume that the error terms, $ \epsilon_t = Y_t - m_t
$, in model (\ref{ModelOfPaper}) form an AR(1) model with possibly
correlated but weak dependent innovations (for precise assumptions
see below), i.e.,
\begin{equation}
\label{ModelOfPaper2}
  \epsilon_t = \rho_t \epsilon_{t-1} + u_t, \qquad t \in \N,
\end{equation}
where $ \rho_t \in (-1,1] $ are unknown parameters. If $ \rho_t =
\rho = 1 $ for all $t$, $ \{ Y_t \} $ is a random walk and {\em
integrated of order $1$, $I(1)$}. Here and throughout the paper we
simply write $ \{ Y_n \} $ if the index set is clear. For $ | \rho | < 1 $ stationary
solutions of the above equation exist.

We consider the following change-point testing problem. The null
hypothesis,
\[
  H_0: \rho_t = 1 \ \mbox{for all $t$},
\]
states that the error terms of the regression model form a random
walk, i.e., are integrated of order $1$. The alternative $ H_1 =
\cup_{q \ge 1} H_1^{(q)} $ with
\[
  H_1^{(q)}: \rho_t = 1, \ t < q, \ \rho_t = \rho, \ t \ge q, \ |\rho|<1
\]
specifies that there exists a change-point such that the subseries
$ \{ \epsilon_t : t \ge q \} $ satisfies stationary AR(1)
equations. It is important to note that the method proposed in
this paper does not require any specification of an alternative.
In Section~\ref{SecCP} we introduce a specific change-point model
related to this testing problem.

Let us consider an example.

\begin{example} Assume $ \epsilon_t = p(L) \xi_t$ with $ \xi_t $ i.i.d. $ N(0,\sigma_\xi^2 )$
for some $ \sigma_\xi > 0 $, $ p(z) = \sum_{j=0}^q \alpha_j z^j $
with coefficients $ \alpha_j \in \R $,  $L$ the lag operator given
by $ L \epsilon_t = \epsilon_{t-1} $. Suppose that the
characteristic polynomial, $ 1 - p(z)$, has exactly one unit root
of multiplicity $1$. Then $ p^*(z) = p(z)/(1-z) $ can be inverted,
and we obtain the representation $ (1-L) \epsilon_t = p^*(L)^{-1}
u_t $, i.e.,
\[
  \epsilon_t = \epsilon_{t-1} + \sum_{j \ge 0} \beta_j u_{t-j},
\]
for certain coefficients $ \beta_j $, see Brockwell and Davis
(1991, Sec. 3.3). Thus, MA($q$)-models with an unit root appear as
a special case for the error terms in model (\ref{ModelOfPaper})
under the null hypothesis.
\end{example}

Concerning the error terms $ \{ u_t \} $ we shall assume the
following mild nonparametric regularity condition making precise
our understanding of weak dependence.

\begin{itemize}
  \item[(E)] $ \{ u_t : t \in \N \} $ is strictly stationary with mean zero
  and $ E | u_1 |^2 < \infty $ such that
  \[
    \sum_{t=1}^\infty | \Cov( u_1, u_{1+t} ) | < \infty,
  \]
  and satisfies a functional central limit theorem (FCLT), i.e.,
  \begin{equation}
  \label{FCLT1}
    T^{-1/2} \sum_{i \le \trunc{Ts}} u_i \stackrel{w}{\to}
    \eta B(s), \qquad T \to \infty,
  \end{equation}
  for some constant $ 0 < \eta < \infty $. Here $B$ denotes a
  Brownian motion with $ B(0) = 0 $, and $ \stackrel{w}{\to} $ stands for weak
  convergence in the Skorohod space $ D[0,1] $. Skorohod spaces
  are briefly discussed at the end of this section.
\end{itemize}

\begin{remark}
\begin{itemize}
\item[(i)] By the Skorohod-Wichura-Dudley representation theorem
(Pollard (1984), Ch. IV.3, Theorem~13), a condition as
(\ref{FCLT1}) is equivalent to the condition: There are Brownian
motions $ B_T $, $ T \ge 1 $, such that
\[
  \sup_{s \in [0,1]} \bigl| T^{-1/2} \sum_{i \le \trunc{Ts}} u_i - \eta
  B_T(s) \bigr| = o_P(1), \qquad T \to \infty.
\]
\item[(ii)] Combining model (\ref{ModelOfPaper2}) with $ \rho_t =
1 $ for all $ t $ under the assumption (E) yields a nonparametric
approach to define the $ I(1)$-property of a time series.
\end{itemize}
\end{remark}

As an example satisfying the assumption (E) let us discuss briefly
$ ARCH(\infty) $ models, a popular parametric class of time series
models.

\begin{example} Recall that $ \{ X_t \} $ satisfies ARCH($\infty$)
equations, if there exists a sequence of i.i.d. non-negative
random variables $ \{ \xi_j \} $, such that
\[
  X_j = \eta_t \xi_t, \qquad \eta_t = a + \sum_{j=1}^\infty b_j
  X_{t-j},
\]
where $ a \ge 0, b_j \ge 0 $ for $ j \in \N $. Suppose now that
\[
  u_t = \sigma_t e_t
\]
where $ \{ e_t \} $ are i.i.d. with $ E(e_t) = 0 $ and $ E(e_t^2)
= 1 $. Put $ \sigma_t = \eta_t $ and $ \xi_t = e_t^2 $ to embed
the classic ARCH model into the above $ ARCH(\infty) $ framework.
Giraitis, Kokoszka and Leipus (2001, Example 2.2 and Theorem 2.1)
have shown that $ \{ u_t \} $ satisfies (E) provided $ E | e_1 |^4
< \infty $ and
\[
  (E \xi_1^4)^{1/4} \sum_{j=1}^\infty b_j < 1.
\]
\end{example}

\subsection{Monitoring Procedure}

\noindent Our stopping time defining the detector essentially
relies on a weighted version of the KPSS test statistic, see
Kwiatkowski et al. (1992), Breitung (2002), and Steland (2007a).
At each time point $ t \le T $ when a new observation is
available, we calculate the full set of residuals $
\widehat{\epsilon}_1(t), \dots, \widehat{\epsilon}_t(t) $ using
all available observations $ Y_1, \dots, Y_t $. Using these
sequentially updated residuals, we calculate an appropriately
weighted version of the unit root test statistic. Define
\[
  U_t = \frac{ t^{-4} \sum_{i=1}^t \bigl( \sum_{j=1}^i
  \widehat{\epsilon}_j(t) \bigr)^2 K( (i-t)/h ) }{ t^{-2}
  \sum_{j=1}^t \widehat{\epsilon}_j^2(t) },
  \qquad t \ge p+1.
\]
In these formulas $ K $, called {\em kernel}, is a nonnegative
function with $ \int K(z) \, dz < \infty $. Kernels such that $
K(z) $ is decreasing for increasing $|z|$ as the Gaussian kernel
or the Epanechnikov kernel given by $ z \mapsto (3/4)
\eins_{[-1,1]}(1-z^2) $, $z \in \R$, have the intuitive appeal
that recent summands get higher weights than past ones. However,
our main results work under the following weak conditions:
\begin{itemize}
  \item[(K1)] $ \| K \|_\infty < \infty $, $ \int K(z) \, dz = 1 $,
  and $ \int zK(z) \, dz = 0 $.
  \item[(K2)] $ K $ is Lipschitz continuous.
\end{itemize}
The parameter $ h = h_T $ is used as a scaling constant in the
kernel and defines the {\em memory} of the procedure. For
instance, if $ K(z) > 0 $ for $ z \in [-1,1] $, and $ K(z) = 0 $
otherwise, $ U_t $ looks back $ h $ observations. We will assume
that
$$
  \lim_{T \to \infty} T/h_T = \zeta \in [ 1, \infty ).
$$
That condition ensures that the number of observations used by the
procedure gets larger as $ T $ increases.

The KPSS or variance ratio control chart is defined as
\[
  R_T = \inf \{ k \le t \le T : U_t \le c_R \}, \qquad  T \ge k,
\]
with the convention $ \inf \emptyset = \infty $. $ T $ is the time
horizon where monitoring stops. For our asymptotic results we
assume $ T \to \infty $, since for applications approximations of
the distribution of $ R_T $  for moderate and large time horizons
$T$ are of interest. $c_R$ is a control limit (critical value)
chosen by the statistician.

It remains to discuss how to choose the control limit $ c_R $.
Since monitoring stops latest at time $T$, we may interpret the
stopping time as a hypothesis test with early stopping in favor of
the alternative. Thus, one may choose $ c_R $ to control
asymptotically the type I error rate of a false decision in favor
of stationarity, i.e.,
\begin{equation}
\label{ChooseCR}
  \lim_{T \to \infty} P_0( R_T \le T ) = \alpha,
\end{equation}
for some given $ \alpha \in (0,1) $. Here $ P_0 $ indicates that
the probability is calculated under the null hypothesis.
Alternatively, one may control a conditional version of the
in-control average run length (CARL). Note that the stopping time
$ R_T $ takes values in the set $ \{ k , \dots, T \} \cup \{
\infty \} $, where $ \infty $ represents no signal, which is the
preferred event under the in-control model. Now we may choose $
c_R $ such that $ \operatorname{CARL}_0 = E_0( R_T | R_T < \infty
) $ is greater or equal to some given value $ \xi \in (k,T) $.
Since our results provide the asymptotic distribution of the
stopping time $ R_T $, one may also choose the control limit to
control other characteristics, e.g., the (conditional) median
average run length. For simplicity of exposition we shall assume
in the sequel that $ c_R $ is chosen such that (\ref{ChooseCR})
holds.

We will assume that monitoring starts after a certain fraction of
the data, i.e.,
\[
  k = \trunc{T \kappa}, \qquad \mbox{for some $ \kappa \in (0,1)
  $},
\]
to avoid that inference is based on too few observations at the
beginning. The event $ R_T \le T $ is interpreted as evidence for
stationary innovations, and we get that information after $ R_T $
observations instead of waiting until time $ T $. If $ R_T =
\infty $, the random walk hypothesis for the error terms is
regarded as compatible with the observed data.

\subsection{Extension to Infinite Time Horizon}

\noindent Suppose we observe sequentially an infinite sequence $
Y_1, Y_2, \dots $ and want to monitor this series with the
detection rule
\[
  \inf \{ k \le t < \infty : U_t \le c_R \}.
\]
In this context, $T$ is simply used to define an appropriate time
scale to determine the bandwidth sequence $ h_T $ with $ T / h_T \to
\zeta $. The FCLT (\ref{FCLT1}) is now replaced by
\[
  \{ T^{-1/2} \sum_{i \le \trunc{Ts}} u_i : s \in [0,\infty) \}
  \stackrel{w}{\to} \{ \eta B(s) : s \in [0,\infty) \},
\]
as $ T \to \infty $, where convergence takes place in the space $
D[0,\infty) $ instead of $ D[0,1] $. All limit theorems in this
paper are formulated for the time interval $ [\kappa,1] $, i.e.,
in the space $ D[\kappa,1] $, but are valid for $ D[\kappa,z] $
for any fixed $ 1 < z < \infty $. $ X_n \stackrel{w}{\to} X $ in $
D[0,\infty) $ is equivalent to
\[
  g_m(t) X_n(t) |_{[0,m]} \stackrel{w}{\to} g_m(t) X(t) |_{[0,m]}
\]
in $ D[0,m] $ for each integer $m$, where $ g_m(t) =
\eins_{[0,m-1]}(t) + (m-t) \eins_{[m-1,m]}(t) $, $ t \in [0,\infty)
$, see Billingsley (1999, Sec. 16) or Pollard (1984, Ch. VI.) Thus,
the results carry over to $ D[\kappa,\infty) $, and there is no loss
in generality to consider the time interval $ [0,\kappa] $.

\subsection{Skorohod Spaces}

\noindent In this paper we will also use the notion of weak
convergence in the Skorohod space $ D([0,1]^2;\R^k)$. Denote the
Skorohod space of cadlag functions $ [0,1] \to \R $ by  $ D[0,1] =
D([0,1];\R) $. Compared to $ D[0,1] $ the space $ D([0,1]^2;\R^k) $
has been only rarely used in the literature. Therefore, we close
this section with a brief exposition of the most important
definitions and facts.

Recall that a sequence $ \{ X, X_n \} $ of random elements with
values in a metric space converges weakly, denoted by $ X_n
\stackrel{w}{\to} X $, as $ n \to \infty $, if $ E h(X_n) \to E
h(X) $, $n \to \infty $, for all measureable real functions which
are bounded and continuous w.r.t. the metric. For a detailed
classic treatment of these issues we refer to Billingsley (1999).

Equip $ D[0,1] $ with the Skorohod metric $d$ yielding a complete
and separable metric space. For $ p\in\N $ let $ D_{\R^p}[\kappa,1] = D(
[\kappa,1]; \R^p ) $ denote the space of all cadlag functions $
[\kappa,1] \to \R^p $ which we equip with the metric $ d_p(f,g) =
\sum_{i=1}^p d(f_i,g_i) $, $ f = (f_1, \dots, f_p)' $, $ g = (g_1,
\dots, g_p)' $, $ f_i, g_i : [0,1] \to \R $, $ i = 1, \dots, p $.
The treatment of its generalization to the index set $ [0,1]^2 $,
i.e., $ D([0,1]^2;\R) $, is more subtle. Let us briefly recall
some facts about this function space and weak convergence of
sequences of $ D( [0,1]^2; \R ) $-valued random elements, as
studied by Straf (1970), Bickel and Wichura (1971), and Neuhaus
(1971). The space $ D( [0,1]^2; \R ) $ can be defined as the
uniform closure of the vector subspace of all simple functions,
i.e., linear combinations of functions of the form $ t \mapsto
1_{E_1 \times E_2}(t) $ where each $ E_i $ is either a
left-closed, right-open subinterval of $ [0,1] $, or the singleton
$ \{ 1 \} $. Here the closure is taken in the space of all bounded
functions $ [0,1]^2 \to \R $.  For functions $ f,g \in D( [0,1]^2;
\R ) $ an appropriate metric, $ d_2(f,g) $, is defined as the
smallest $ \varepsilon > 0 $ such that there exist continuous
bijections $ \lambda_1, \lambda_2 : [0,1] \to [0,1] $ with $ \|
\lambda - \id \|_\infty \le \varepsilon $ and $ \| f - g \circ
\lambda \|_\infty \le \varepsilon $. Here $ \lambda = (\lambda_1,
\lambda_2) $ and $ g \circ \lambda(r,s) = g(\lambda_1(r),
\lambda_2(s)) $ for $ (r,s) \in [0,1]^2. $ A sequence $ \{ f_n \}
\subset D([0,1]^2;\R) $ converges to some $ f \in D([0,1]^2;\R) $
iff there exists some sequence $ \{ \lambda_n \} $ of pairs of
continuous bijections $[0,1] \to [0,1] $ such that $ \| f_n \circ
\lambda_n - f \|_\infty \to 0 $ and $ \| \lambda - \id \|_\infty
\to 0 $, as $ n \to \infty $. Further, if $ f \in C([0,1]^2;\R) $
convergence in the Skorohod metric implies uniform convergence,
since in this case $f$ is uniformly continuous. It turns out that
$ (D([0,1]^2;\R),d_2) $ is a separable metric space, a common
framework to define weak convergence of probability measures and
random elements.

\section{ASYMPTOTIC RESULTS FOR INTEGRATED \\ PROCESSES}
\label{SecIP}

\noindent This section is devoted to a detailed study of the
proposed procedure unter the null hypothesis that the error terms
of the regression model behave as a random walk. Our approach is
to represent the KPSS control chart as an inf-functional of the
stochastic process associated to the sequence $ \{ U_t \} $. That
process turns out to be a functional of the stochastic process
associated to the residuals up to negligible terms. We provide
functional central limit theorems for these processes and a
central limit theorem for the stopping time $ R_T $.

We need some notations. Let $ \matX_n $ denote the design matrix
for a polynomial regression of order $p$ with intercept based on
$n$ observations, i.e.,
\[
  \matX_n = \left[
   \begin{array}{cccc}
     1      &  1      &  \cdots  &  1 \\
     1      &  2      &  \cdots  &  2^p \\
     \vdots &  \vdots &          &  \vdots \\
     1      &   n     &  \cdots  &  n^p
   \end{array}
   \right]
   = [ \vecx_1, \dots, \vecx_n ]',
\]
where
\[
  \vecx_t = (1, t, \dots, t^p )'.
\]
 Define for $ p+1 \le t \le T $ the random vectors
\begin{eqnarray*}
  \bfeps_t & = & ( \epsilon_1, \dots, \epsilon_t )' \\
  \widehat{\bfeps}_t
           & = & ( \matid_t - \matX_t( \matX_t' \matX_t )^{-1} \matX_t' ) \vecY_t,
\end{eqnarray*}
with $ \vecY_t = (Y_1, \dots, Y_t)' $. $ \matid_t $ denotes the
$t$-dimensional identity matrix.

\subsection{Residual Process without Updating}

\noindent Let us first consider the natural process associated to
the sequence $ \widehat{\epsilon}_1, \dots, \widehat{\epsilon}_T $
of residuals, where at time $t$ the current residual $
\widehat{\epsilon}_t $ is simply added to the residuals $
\widehat{\epsilon}_i $, $  \trunc{ T\kappa } \le i < t $. Here the
former residuals are {\em not} updated. In the sequel $ \trunc{Ts}
$ stands for the current time point. The stochastic process
associated to $ \widehat{\epsilon}_1, \dots, \widehat{\epsilon}_T
$ is given by
\[
  \widehat{E}_T(s) = T^{-1/2} \widehat{\epsilon}_{ \trunc{Ts} }, \qquad s \in [ \kappa, 1 ],
\]
where $ \widehat{\epsilon}_t = 0 $ for $ 0 \le t < p+1 $, and
\[
  \widehat{\epsilon}_{\trunc{Ts}} = Y_{\trunc{Ts}} - \vecx_{\trunc{Ts}}'
   ( \matX_{\trunc{Ts}}' \matX_{\trunc{Ts}} )^{-1} \matX_{\trunc{Ts}}' \vecY_{\trunc{Ts}}
\]
is the last coordinate of the vector $
\widehat{\bfeps}_{\trunc{Ts}} = ( \widehat{\epsilon}_1, \dots,
\widehat{\epsilon}_{\trunc{Ts}} )' $.

We have to introduce the weighting matrix
\[
  \matW_t = \diag( 1, t^{-1} ,\dots, t^{-p} ), \qquad t \in \N,
\]
to take into account the order of the polynomial regressors.

\begin{lemma}
\label{Lemma1} Fix $ \kappa \in (0,1) $. Assume (E). Then
  \[
    T^{-3/2} \matW_\trunc{Ts} \matX_{\trunc{Ts}}' \bfeps_{\trunc{Ts}}
    \stackrel{w}{\to} \eta \int_0^s (1,r/s,\dots,(r/s)^p)' B(r) \, dr,
    \qquad \mbox{in $D_{\R^p}[\kappa,1] $},
  \]
  as $ T \to \infty $, where the limit is almost surely (a.s.) continuous.
\end{lemma}

Lemma~\ref{Lemma1} plays a crucial role in the proofs of the main
results, but it is also interesting in its own right. Notice that $
\matX_{\trunc{Ts}}' \bfeps_{\trunc{Ts}} $ is the natural
sufficient statistic when the errors are i.i.d. normal. The lemma
states that for random walk error terms with weak dependent
increments  the correct scaling operator for the natural
sufficient statistic is given by $ T^{-3/2} \matW_{\trunc{T\cdot}} $ to obtain a
non-degenerate distributional limit. The limit process is given by
the vector of weighted integrals of Brownian motion, $ \eta
\int_0^s (r/s)^k B(r) \, dr $, $ k = 0, \dots, p $, where the
integral is a Riemann integral. The factor $ \eta $ summarizes the
impact of the correlation of the increments.

Let us introduce the Hilbert matrix of dimension $ p+1 $ given by
\[
  \matH = ( 1/(i+j-1) )_{i,j \in \{ 1, \dots,p+1 \}}.
\] It is known that its inverse, $ \matH^{-1} $,  has entries
\[
  ( \matH^{-1} )_{i,j} = (-1)^{i+j} (i+j-1) { p+i \choose p+1-j }{ p+j \choose p+1-i }
  { i+j-2 \choose i-1 }^2,
\]
see Choi (1983).

We need the following simple result about sufficient conditions
for uniform convergence of the inverse of a sequence of invertible
matrix-valued functions $ \matA_n(x) $, $ \matA_n : \R \to \R^{l
\times l} $,  to the inverse of its limit $ \matA(x) $. Let $ \|
\cdot \| $ denote the Euclidean vector and matrix norm,
respectively.

\begin{lemma}
\label{Lemma2}
 Suppose $ \{ \matA(x), \matA_n(x) : n \ge 1 \} $,
is a sequence of $k$-dimensional matrix-valued functions such that
\[
  \sup_x \| \matA_n(x) - \matA(x) \|_2 = o(1).
\]
If
\begin{equation}
\label{VsM}
  0 < \inf_x \sigma_1(x) \quad \mbox{and} \quad \sup_x \sigma_k(x) <
  \infty,
\end{equation}
where $ \sigma_1(x) $ ($ \sigma_k(x) $) denotes the smallest
(largest) eigenvalue of $ \matA(x)^* \matA(x) $, then
\[
  \sup_x \| \matA_n^{-1}(x) - \matA^{-1}(x) \|_2 = o(1).
\]
\end{lemma}

\begin{theorem}
\label{SeqRes1} Fix $ \kappa \in (0,1) $. Assume (E). Then, under
the null hypothesis $ H_0 $,
  $$
    \widehat{E}_T \stackrel{w}{\to} \calE,
    \qquad \mbox{in $ D[\kappa,1] $},
  $$
  as $ T \to \infty $, where the a.s. continuous process $ \calE $ is given by
  $$
    \calE(s) = \eta \left\{ B(s) -
      s^{-1} \eins' \matH^{-1} \int_0^s ( 1, r/s, \dots, (r/s)^p )' B(r) \, dr \right\},
  $$
  for $ s \in [\kappa,1] $.
\end{theorem}

This theorem provides an explicit formula for the limit process of
$ \widehat{E}_T $. The limit process is a linear function of
Brownian motion $ B(s) $ and the limit process appearing in
Lemma~\ref{Lemma1}.

\subsection{Sequentially updated Residual Process}

\noindent Again, $ \trunc{Ts} $ denotes the current time point and
$ \trunc{Tr} $ stands for another time point, in most cases a
previous one. Let us now consider the two-parameter stochastic
process
$$
  \widehat{E}_{\trunc{Tr}}( \trunc{Ts} )
  = T^{-1/2} \widehat{\epsilon}_{\trunc{Tr}}( \trunc{Ts} ),
  \qquad r \in [\kappa,s], \ s \in [\kappa,1],
$$
where for $ p+1 \le k \le t \le T  $ we denote by $
\widehat{\epsilon}_k(t) $ the $k$-th residual associated to the
observation $ Y_k $, calculated using the data $ Y_1, \dots, Y_t
$. This means, having observed the $n$th observation, {\em all}
residuals are updated. We call $ \widehat{E}_{\trunc{Tr}}(
\trunc{Ts} ) $ the {\em sequentially updated residual process}.
Extend the definition by putting $ \widehat{E}_{\trunc{Tr}}(
\trunc{Ts} ) = 0 $ if $ r > s $ or $ r,s \in [0,\kappa) $ to
obtain a $ D([0,1]^2)$-valued process. Note that
\[
  \widehat{\epsilon}_{\trunc{Tr}}( \trunc{Ts} )
  = Y_{\trunc{Tr}} - \vecx_{\trunc{Tr}}' \bhbeta_{\trunc{Ts}},
\]
with
\begin{align}
\label{DefBetaHat}
  \bhbeta_{\trunc{Ts}} &= (\matX_{\trunc{Ts}}'\matX_{\trunc{Ts}})^{-1} \matX_{\trunc{Ts}}'\bfeps_{\trunc{Ts}}
  \\ \nonumber
  \bfeps_{\trunc{Ts}}  &= (\epsilon_1,\dots,\epsilon_{\trunc{Ts}})'.
\end{align}

\begin{theorem}
\label{FCLTUpdatedResiduals} Fix $ \kappa \in (0,1) $ and assume
(E). We have under $ H_0 $
$$
  \widehat{E}_{\trunc{Tr}}( \trunc{Ts} )
  \stackrel{w}{\to} \calE(r,s),
$$
in $ D([\kappa,1]^2;\R) $, as $ T \to \infty $, where the process $
\calE $ is given by
\begin{equation}
\label{DefE}
  \calE(r,s) = \eta \left\{
    B(r) - \vecv(r,s) s^{-1} \matH^{-1} \int_0^s \vecv(u,s) B(u) \,
    du
  \right\}
\end{equation}
with
\begin{equation}
  \vecv(r,s) = ( 1, r/s, \dots, (r/s)^p)',
\end{equation}
for $ \kappa \le r \le s \le 1 $.
\end{theorem}

Notice that the limit process for sequentially updated
residuals has a similar structure as for residuals without
updating, but the vector functions appearing in the definition of
$ \calE(r,s) $ now depend on both $r$ and $s$. Again, the
influence of the dependence structure of the error terms is
summarized by the factor $\eta $.

\begin{example} Explicit representations of the limit processes are now easy
to obtain. Let us consider dimensions $p=1$ and $p=2$, which are
of special importance for applications.
\begin{itemize}
\item[(i)] For $ p = 1 $ we have $ \matH^{-1} = \left[
\begin{array}{cc} 4 & -6 \\ -6 & 12 \end{array} \right] $ and
$$
  s^{-1} \matH^{-1} \left[ \int_0^s B(r) \, dr, \int_0^s r B(r) \, dr \right]' =
  \left(
  \begin{array}{c}
    \frac{4}{s} \int_0^s B(r) \, dr - \frac{6}{s^2} \int_0^s r B(r) \, dr \\
   -\frac{6}{s} \int_0^s B(r) \, dr + \frac{12}{s^2} \int_0^s r B(r) \, dr
  \end{array}
  \right).
$$
Thus,
$$
    \mathcal{E}(r,s) = \eta \left\{ B(r) +
     \left( \frac{6r}{s^2} - \frac{4}{s} \right) \int_0^s B(u) \, du
     + \left( \frac{6}{s^2} - \frac{12 r}{s^3} \right) \int_0^s u B(u) \, du \right\},
$$
for $ r,s \in [\kappa,1]$, $ r \le s $.
\item[(ii)] If $ p = 2 $,
$ \matH^{-1} \int_0^s [1,r/s,r^2/s^2]' B(r) \,dr $ is given by
$$
  \left(
    \begin{array}{c}
    9 \int_0^s B(r) \, dr - \frac{36}{s} \int_0^s r B(r) \, dr + \frac{30}{s^2} \int_0^s r^2 B(r) \, dr \\
    -36 \int_0^s B(r) \, dr + \frac{192}{s} \int_0^s r B(r) \, dr - \frac{180}{s^2} \int_0^s r^2 B(r) \, dr \\
    30 \int_0^s B(r) \, dr - \frac{180}{s} \int_0^s r B(r) \, dr + \frac{180}{s^2} \int_0^s r^2 B(r) \, dr
    \end{array}
  \right)
$$
We obtain
\begin{eqnarray*}
\mathcal{E}(r,s) &=& \eta \biggl\{ B(s) -
 \biggl( \frac{9}{s} - \frac{36 r}{s^2} + \frac{30r^2}{s^3} \biggr) \int_0^s B(r) \, dr  \\
&& \quad - \biggl( -\frac{36}{s^2} + \frac{192r}{s^3} - \frac{180 r^2}{s^4} \biggr) \int_0^s r B(r) \, dr \\
&& \quad - \biggl( \frac{30}{s^3} - \frac{180r}{s^4} +
\frac{180r^2}{s^5} \biggr) \int_0^s r^2 B(r) \, dr \biggr\},
\end{eqnarray*}
for $ r,s \in [\kappa,1] $, $ r \le s $.
\end{itemize}
\end{example}

\begin{remark} Based on these explicit formulas, simulating trajectories of
the process $ \calE(r,s) $ becomes a feasible task. Using
Donsker's theorem one may simulate trajectories of $B(r)$ and
employ numerical integration to simulate the moment functions $
\int_0^s r^k B(r)\,dr $, $ s \in [\kappa,1] $, $k \in \N $, appearing in
the formulae.
\end{remark}

\subsection{Weighted Variance Ratio Process}

\noindent We are now in a position to examine the process
associated to the sequence of control statistics $ \{ U_t \} $.
For brevity of exposition we present the results for the
sequentially updated residuals. The required modifications when
using the sequential residuals without updating are
straightforward.

Define the {\em kernel-weighted variance ratio process}
\[
  V_T(s) = \frac{ \trunc{Ts}^{-4} \sum_{i=g_T}^{\trunc{Ts}} \bigl( \sum_{j=1}^i
  \widehat{\epsilon}_j( \trunc{Ts} ) \bigr)^2 K( (i - \trunc{Ts})/h ) }
  { \trunc{Ts}^{-2}
  \sum_{j=g_T}^{\trunc{Ts}} \widehat{\epsilon}_j^2( \trunc{Ts} ) },
  \qquad s \in [0,1].
\]
Here and in the sequel we agree to put $0/0=0$. $ g_T $ denotes
the time point where calculations start. To ensure both that the
residuals can be calculated and the sums appearing in the
definition of $ V_T(s) $ have a reasonable number of summands for
all $ s \in [\kappa,1] $, we assume  $ p + 1 \le g_T <
\trunc{T\kappa} $. A plausible choice is
\[
  g_T = \trunc{T \gamma}, \qquad \mbox{for some $ \gamma \in (0,\kappa) $}.
\]
Then $ g_T/T \to \gamma > 0 $. More generally, let
\begin{equation}
\label{GammaLimit}
  \gamma = \lim_{T \to \infty} g_T/T \in [0,\kappa].
\end{equation}
The stopping time $ R_T $ can now be represented as
\[
  R_T = T \inf \{ s \in [\kappa,1] : V_T( s ) \leq c \}.
\]

We are now in a position to formulate the main result.

\begin{theorem}
\label{KPSS}
 Fix $ \kappa \in (0,1) $ and assume (E). Under $ H_0 $ we have
\[
  V_T(s) \stackrel{w}{\to} \calV(s) = \frac{ s^{-2} \int_\gamma^s
  \bigl( \int_0^r \calE(s,t) \, dt \bigr)^2 \, K(\zeta(r-s)) dr }
  { \int_\gamma^s \calE^2(s,r) \, dr }, \qquad T \to \infty,
\]
in the space $ D[\kappa,1] $. The limit process is continuous w.p.
$1$ and depends only on $K$, $ \zeta $, and Brownian motion $B$, but
not on the quantity $ \eta $.
\end{theorem}

We discuss this theorem at the end of this section in greater
detail.

\subsection{KPSS (Variance Ratio) Residual Control Chart}

\noindent The central limit theorem for the stopping time $R_T $
of the KPSS residual control chart appears now as a corollary to
Theorem~\ref{KPSS}.

\begin{corollary}
\label{TheCorollary}
  For the stopping time $R_T$ we have under the conditions of
  Theorem~\ref{KPSS}
  \[
    R_T/T \stackrel{d}{\to} \calR = \inf \{ s \in [\kappa,1] : \calV(s) \le c_R \},
  \]
  as $ T \to \infty $.
\end{corollary}

As a consequence, the KPSS residual control chart can be designed
to achieve a given nominal significance level $ \alpha \in (0,1)
$. Indeed, Corollary \ref{TheCorollary} implies that $ P_0( R_T \le T ) \to P_0(
\calR \le 1 ) $, as $ T \to \infty $. Since
\[
  \calR \le 1 \Leftrightarrow \inf_{s \in [\kappa,1]} \calV(s) \le
  c_R,
\]
we select the control limit as
\[
  c_R = F^{-1}(1-\alpha),
\]
where $ F $ denotes the distribution function of $ \inf_{s \in
[\kappa,1]} \calV(s) $.

\begin{remark} Having in mind practical applications it is worth
discussing the following issues.
\begin{itemize}
  \item[(i)] $ V_T $ converges weakly to a stochastic process which
does not depend on any nuisance parameter. When a kernel $K$ and
the parameter $ \zeta $ are selected, the process $ \calV $ is
known. This means, the asymptotic law of $ V_T $ is
distribution-free. As a consequence, the asymptotic distribution
of $ R_T $ is also asymptotically distribution-free. \item[(ii)]
In practice, one can simulate trajectories from the limit process
and calculate for each trajectory the time point where the control
limit $ c_R $ is reached. In this way one can simulate the
asymptotic distribution of $ R_T $ to determine a control limit $
c_R $ such that the resulting asymptotic type I error rate is $
\alpha $.
\end{itemize}
\end{remark}

\section{ASYMPTOTIC RESULTS FOR A CHANGE-POINT \\ MODEL}
\label{SecCP}

\noindent The results of the previous section allow to design
monitoring procedures and to study the behavior of the resulting
procedure under the null hypothesis (in-control model) that the
underlying time series of observations follows a polynomial
regression model with random walk error terms under the stated
regularity assumptions.

In this section we discuss the asymptotic behavior of the KPSS
residual monitoring approach under a change-point model, where the
first part of the time series behaves as a random walk and the
second part is stationary. We assume
\begin{equation}
\label{CPM}
  \epsilon_t = \left\{
  \begin{array}{cc}
    \sum_{j=0}^t u_j, & t = 0, \dots, \trunc{T \vartheta} -
    1, \\
    \xi_T u_t, & t = \trunc{T \vartheta}, \dots, T.
  \end{array}
  \right.
\end{equation}
After the change-point $ q = \trunc{T \vartheta} $, which is given
by the fixed but unknown parameter $ \vartheta \in (0,1) $, the
error terms change and are no longer a random walk. $ \{ \xi_T \}
$ is a  sequence of scale constants  satisfying the condition
\begin{equation}
\label{CondXi}
  \xi_T \sim T^{\beta}, \qquad \mbox{for some $ \beta \ge 0 $}.
\end{equation}
We shall need a further constraint on $ \beta $ which will be
discussed below. If $ \beta = 0 $, the error process after the
change, i.e., $ \{ \epsilon_t : q \le t \le T \} $ with $q =
\trunc{T\vartheta} $, is stationary. However, we allow for positive
values of $ \beta $. In this case the error terms form a row-wise
stationary array. For simplicity of exposition, we omit the
dependence of $ \epsilon_t $ on $T$ in our notation.

Our asymptotic results require the following additional assumptions.
\begin{itemize}
  \item[(C1)] $ \{ u_t \} $ is a strictly stationary process with
  \[
    \lim_{x \to \infty} \frac{ P( | u_1 | > x ) }{ x^{-\gamma} } <
    \infty,
  \]
  for some $ \gamma > 2 $ and satisfies the FCLT
  \begin{equation}
    T^{-1/2} \sum_{i \le \trunc{Ts}} u_i \stackrel{w}{\to}
    \eta B(s), \qquad T \to \infty,
  \end{equation}
  for some constant $ 0 < \eta < \infty $, where again $B$ denotes
  standard Brownian motion starting at $0$.
  \item[(C2)] The parameters $ \alpha $ and $ \beta $ satisfy the relations $ 0 \le \beta < 1/2 $ and $ \gamma >
  \frac{1}{1/2 - \beta} $.
\end{itemize}

Note that the condition on the tail probabilities ensures that the
$ E|u_t|^2 < \infty $.

In the sequel, we use the same notation for the quantities defined
for the polynomial regression model with error terms $ \{ \epsilon_t
\} $ satisfying the change-point model above.

Let us again start with the residual process. We only discuss the
FCLT for the process of sequentially updated residuals, $
\widehat{E}_{\trunc{Tr}}( \trunc{Ts} ) $, $ \kappa \le r \le s \le
1 $, which is defined as before.

\begin{theorem}
\label{FCLTUpdatedResidualsCP} Suppose the change-point model
(\ref{CPM}) holds. Additionally, assume that (E), (C1), and (C2)
are satisfied. Then, for any fixed $ \kappa \in (0,1) $, the
following assertions hold true.
\begin{itemize}
\item[(i)] We have in the space $ D([\kappa,1];\R^p) $,
\[
  T^{-3/2} \matW_{\trunc{Ts}} \matX_{\trunc{Ts}}'
  \bfeps_{\trunc{Ts}}
  \stackrel{w}{\to}
  \eta \int_0^{s} (1,r/s, \dots, (r/s)^p )' B(r)
  \, dr \eins_{[\kappa,\vartheta)}(s),
\]
as $ T \to \infty $. \item[(ii)] The sequentially updated LS
residual process converges weakly,
  $$
    \widehat{E}_T \stackrel{w}{\to} \calE_{\vartheta},
    \qquad \mbox{in $ D([\kappa,1]^2;\R) $},
  $$
  as $ T \to \infty $, where the cadlag process $ \calE $ is given by
  $$
    \calE_{\vartheta}(r,s) = \eta \left\{ B(s) -
      s^{-1} \vecv(r,s) \matH^{-1} \int_0^s \vecv(u,s) B(u) \, du
      \right\} \eins_{[\kappa,\vartheta)}(s),
  $$
  for $ \kappa \le r \le s \le 1 $.
\end{itemize}
\end{theorem}

The next result shows that under the aforementioned conditions the
asymptotic distribution of the  kernel-weighted variance ratio
process is obtained by replacing formally $ \calE $ by $
\calE_\vartheta $ in the limit process.

\begin{theorem}
\label{KPSS_CP} \label{FCLTCP} Suppose the change-point model
(\ref{CPM}), assumption (E), (C1) and (C2) are satisfied. Then,
for any fixed $ \kappa \in (0,1) $,
\[
   V_T(s) \stackrel{w}{\to} \calV_\vartheta(s) = \frac{ s^{-2} \int_\gamma^s
     \bigl( \int_0^r \calE_\vartheta(s,t) \, dt \bigr)^2 \, K(\zeta(r-s)) dr }
     { \int_\gamma^s \calE^2_\vartheta(s,r) \, dr }, \qquad T \to \infty,
\]
in the space $ D[\kappa,1] $. The limit process depends only on
$K$, $ \zeta $, and Brownian motion $B$, and the change-point
parameter $ \vartheta $.
\end{theorem}

Again, the central limit theorem for the KPSS residual control
chart under the change-point model appears as a corollary.

\begin{corollary} Under the assumptions of Theorem~\ref{KPSS_CP}, the stopping time $R_T$
satisfies
\[
  R_T/T \stackrel{d}{\to} \calR_\vartheta = \inf \{ s \in [\kappa,1] : \calV_\vartheta(s) > c_R \},
\]
as $ T \to \infty $.
\end{corollary}

\section{SIMULATIONS}
\label{SecSim}

\noindent We conducted a Monte Carlo study to investigate the
properties of the KPSS monitoring procedure when applied to
residuals. Time series of length $ T = 500 $ according to model
\[
  Y_t = \beta_0 + \beta_1 \cdot t + (\beta_1 + \Delta) t
  \eins_{\{q,q+1,\dots\}}(t) + \epsilon_t,
\]
where
\[
  \epsilon_t = \left\{
    \begin{array}{ll}
      \sum_{i=1}^t \eta_i, \qquad & \text{for $ t < q$}, \\
      \sum_{i=1}^{q-1} \eta_i + \eta_t, \qquad & \text{for $ q \le
      t \le T $},
    \end{array}
  \right.
\]
with
\[
  \eta_t = \rho \eta_{t-1} + \xi_t - \beta \xi_{t-1}, \qquad \xi_t
  \stackrel{i.i.d.}{\sim} N(0,1),
\]
were simulated. Let us first discuss the construction of the
innovation terms $ \eta_t $. The AR parameter was chosen as $ \rho
= 0.3 $ and the MA parameter $ \beta $ from the set $ \{ -0.8, 0,
0.8 \} $. Thus, $ \{ \eta_t \} $ is a correlated but weakly
dependent sequence with mean $0$. For time points $ t < q $ the
obervations $ Y_t $ are given by a random walk with correlated
increments $ \eta_i $. At the change-point $ q $ the process
changes its behavior. The random walk stops and correlated error
terms $ \eta_t $ determine the behavior of $ \epsilon_t $.

Concerning the design of the monitoring procedures we used the
Gaussian kernel, $ K(z) = (2 \pi)^{-1} \exp(-z^2/2) $, $ z \in \R
$, and the bandwidths $ h \in \{ 25, 50 \} $, yielding $ \zeta \in
\{ 20, 10 \} $. The deterministic component of the model is given
by a linear trend whose slope, depending on $ \Delta $, may change
at the change-point, too. If $ \Delta \not= 0 $, there is both a
change in the error terms and a change in the slope. That should
make the detection of the change to stationarity of the errors
more difficult, since the residuals are estimated assuming a
constant slope.

In a first step, we examined for the setting $ \Delta = 0 $, $ h =
25 $, and $ T = 500 $, the relationship between the control limit $
c $ and, firstly, the probability that the method gives a signal
(Figure 1) and, secondly, the conditional average run length (CARL)
given that the method gives a signal at all (Figure 2). Since
monitoring stops latest at the $500$th observation, trajectories
crossing the control limit later are not taken into account. The
CARL is the average run length corresponding to all trajectories
yielding a signal until time $ 500 $. The curves, which can also be
used to choose the control limit, are quite similar for $ \beta \in
\{ -0.8, 0 \} $, but there is an effect for positive values of $
\beta $. For the considered setting it also becomes apparent that
common type I error rates correspond to rather large CARL values. On
the other side, if the procedure is designed to yield CARL values
of, say, $ 300 $, the chart works on a type I error rate which is
usually regarded as unacceptable from a hypothesis testing
viewpoint. However, note that this is partly due to the fact that we
studied monitoring with a time horizon. Without a time horizon the
average run lengths would be substantially higher yielding smaller
control limits and, as a consequence, smaller associated type I
error rates.

\begin{figure}\label{Plot1}
\centering
\includegraphics[width=6cm]{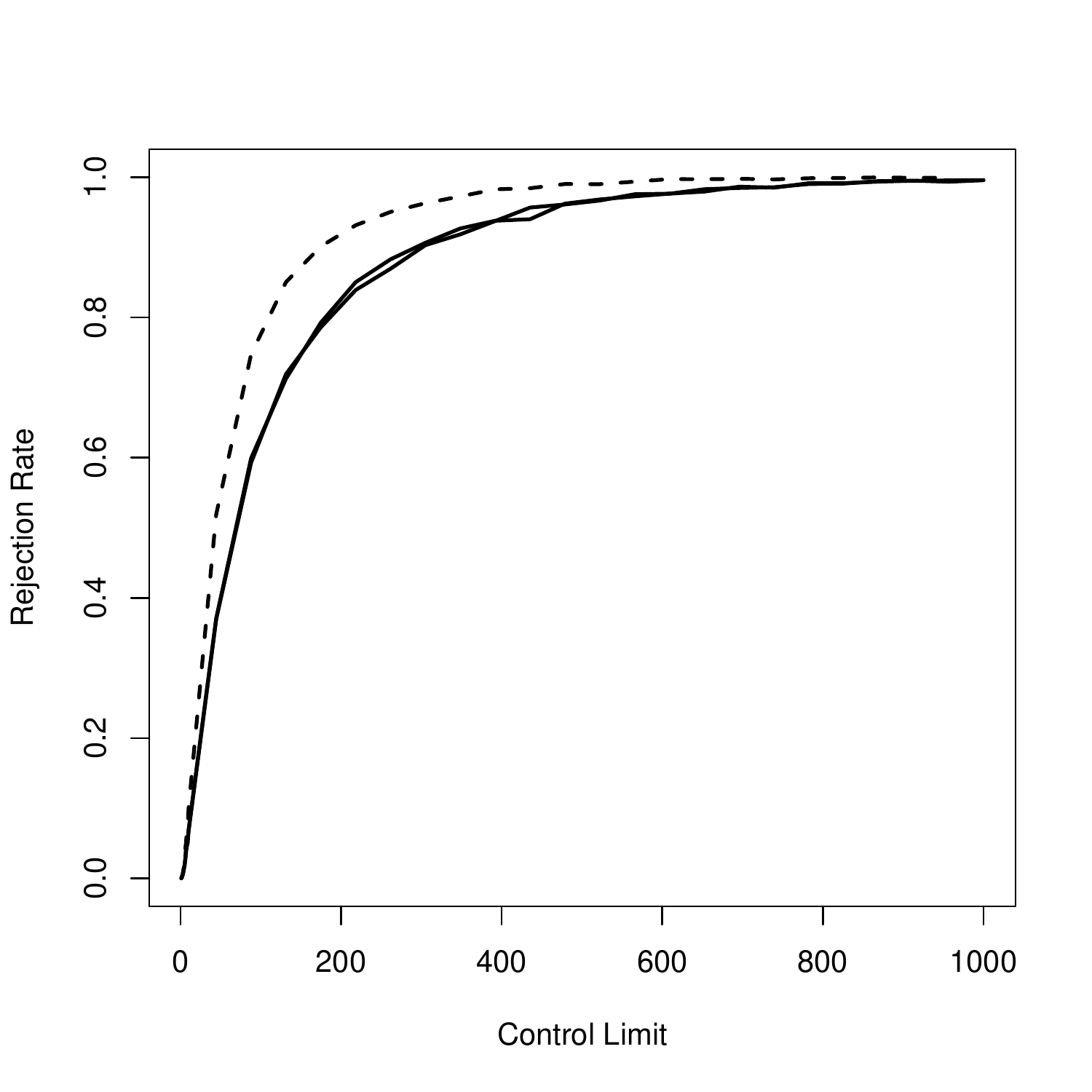}

\textit{\textbf{Figure~1. }}{Empirical rejection rates as a function of $ 10^6 $ times the control limit $c$
for $ \beta = -0.8 $, $ \beta = 0 $, and $ \beta = 0.8 $ (dashed).}
\end{figure}

\begin{figure}\label{Plot2}
\centering
\includegraphics[width=6cm]{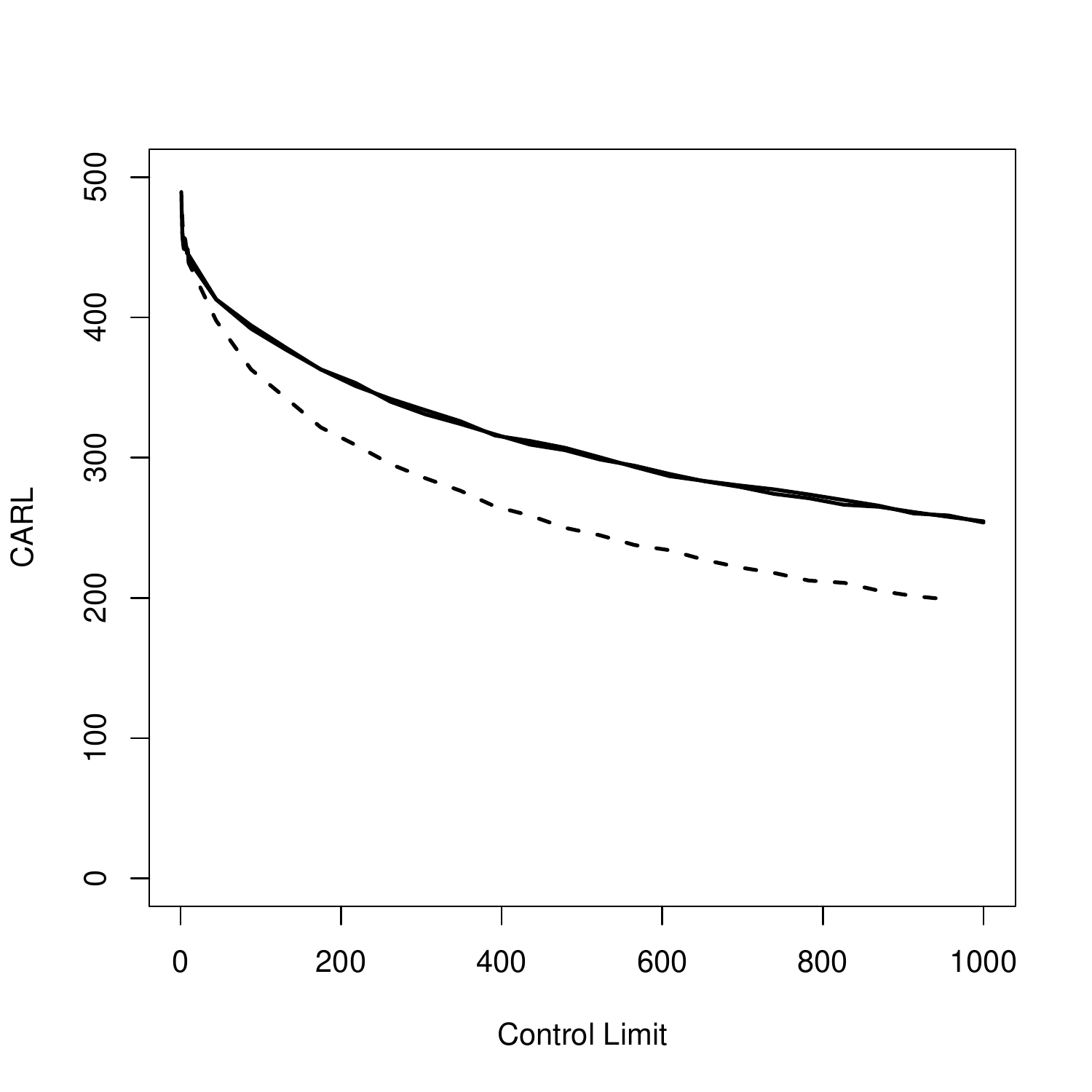}

\textit{\textbf{Figure~2. }}{Conditional average run length (CARL) as a function of $
10^6 $ times the control limit $c$ for $ \beta = -0.8 $, $ \beta =
0 $, and $ \beta = 0.8 $ (dashed).}
\end{figure}

We also simulated the power of the KPSS variance ratio residual
control chart when designed to achieve a type I error rate of $
\alpha = 0.05 $. The corresponding control limit was obtained by
simulating from the limit distributions. We examine the cases $
\Delta = 0 $ and $ \Delta = 0.5 $, where the latter case
corresponds to a change to stationary errors term with an
additional change of the slope.

Table~1 provides the simulated rejection rates. It can be seen
that the KPSS control chart is quite robust with respect to the
parameter $ \beta $ determining the degree of correlation for the
increments of the random walk. That behavior is consistent with
the theoretical and empirical findings in Steland (2007a), where
related monitoring procedures for time series without trends have
been studied in detail. The results show that an early change can
be detected quite well, but late changes are very hard to detect.
However, this is, of course, a problem for all methods, and for
the statistical problem at hand, $ N = 500 $ is not a large
maximal number of observations. The results also indicate that the
power is quite robust with respect to additional changes in slope.

We may summarize that the KPSS control chart for residuals
provides a quite reliable tool to detect stationary errors in
polynomial regression models.

\begin{center}
\textbf{Table~1.} {Empirical rejection rates of the KPSS control
chart. The right columns ($\Delta = 0.25$) correspond to a change
of the slope}

\begin{tabular}{ccccccc} \hline
              &  \multicolumn{3}{c}{$ \Delta = 0 $}       &
              \multicolumn{3}{c}{$ \Delta = 0.25 $}  \\ \hline
              &   \multicolumn{3}{c}{$ \beta $}  &
              \multicolumn{3}{c}{$\beta$} \\
change-point  & $-0.8$ & $0$ & $ 0.8$ & $ -0.8 $ & $0$ & $ 0.8 $
              \\ \hline
Results for  $ h = 25 $ \\
  $ 25 $  & $ 0.44 $  & $ 0.60 $ & $0.94$ & $ 0.45 $ & $ 0.57 $ & $ 0.94 $ \\
  $ 75 $  &  $0.16$   & $ 0.18$  & $0.41$ & $ 0.18 $ & $ 0.19 $ & $ 0.38 $ \\
  $100$   &  $0.17$   & $0.16$   & $0.29$ & $ 0.16 $ & $ 0.16 $ & $ 0.29 $ \\
no-change &  $0.06$   & $0.06$   & $0.07$ & $ 0.06 $ & $ 0.06 $ & $ 0.10 $ \\
\hline

Results for $ h = 50 $ \\\
  $ 25 $  & $ 0.52 $  & $ 0.61 $ & $0.97$ & $ 0.53 $ & $ 0.60 $ & $ 0.97 $ \\
  $ 75 $  &  $0.18$   & $ 0.19$  & $0.44$ & $ 0.17 $ & $ 0.18 $ & $ 0.44 $ \\
  $100$   &  $0.15$   & $0.16$   & $0.30$ & $ 0.15 $ & $ 0.16 $ & $ 0.33 $ \\
no-change &  $0.03$   & $0.03$   & $0.06$ & $ 0.03 $ & $ 0.03 $ & $ 0.07 $ \\
\hline

\end{tabular}
\end{center}

\section*{ACKNOWLEDGMENTS}

\noindent The author gratefully acknowledges a grant (STE
1034/6-1) from Deutsche Forschungsgemeinschaft (DFG) and thanks
the Editor, an Associate Editor, and a referee for the excellent
review process and valuable remarks. He also thanks S. Teller for
her careful proofreading, H. Satvat and H. Schottm\"uller for
building up a Linux cluster for parallel computing which was used
for simulations, and B. Giese for preparing the final Latex
layout.

\section*{APPENDIX A: PROOFS OF RESULTS FROM SECTION 3}

\begin{myproofa}{of Lemma~\ref{Lemma1}}{1}
  Note that for each $ s \in [\kappa,1] $ we have
  $$
    \matX_{\trunc{Ts}}' \bfeps_{\trunc{Ts}}
    = \left( \sum_{t=1}^{\trunc{Ts}} \epsilon_t, \sum_{t=1}^{\trunc{Ts}} t \epsilon_t,
     \dots, \sum_{t=1}^{\trunc{Ts}} t^p \epsilon_t \right)'
  $$
  yielding
  \begin{align}
  \label{RepMatr}
    \nonumber
    T^{-3/2} \matW_{\trunc{Ts}} \matX_{\trunc{Ts}}' \bfeps_{\trunc{Ts}}
    &=
    T^{-3/2} \left( \sum_{t=1}^{\trunc{Ts}} (t/\trunc{Ts})^{i-1} \epsilon_t \right)_{i=1,\dots,p+1} \\
    &= T^{-1/2} \left( \int_0^s (\trunc{Tr}/\trunc{Ts})^{i-1} \epsilon_{\trunc{Tr}} \, dr \right)_{i=1,\dots,p+1} \tag{A.1}\\
    &= \left( \int_0^s (\trunc{Tr}/\trunc{Ts})^{i-1} T^{-1/2} \epsilon_{\trunc{Tr}} \, dr
    \right)_{i=1,\dots,p+1}. \nonumber
  \end{align}
  It is straightforward to check that
  \begin{equation*}\tag{A.2}
  \label{PowerFrac}
    \sup_{\kappa \le r \le s \le 1} \max_{1 \le i \le p} | ( \trunc{Tr} / \trunc{Ts} )^i - (r/s)^i | =
    O(1/T).
  \end{equation*}
  Hence,
  \[
   \sup_{0 \le r \le s \le 1 } \max_{1 \le i \le p} \bigl| \int_0^s (\trunc{Tr}/\trunc{Ts})^{i-1} z(r) \, dr -
   \int_0^s (r/s)^{i-1} z(r) \, dr
   \bigr| = O(1/T).
  \]
  If we define the functional $ \tau : (D[\kappa,1],d) \to (D_{\R^p}[\kappa,1], d_p) $
  by
  \[
      \tau( z )( s ) =
    \biggl(
      \int_0^s (r/s)^{i-1} z(r) \, dr \biggr)_{i=1,\dots,p+1},
    \qquad s \in [\kappa, 1],
  \]
  for any $ z \in D[\kappa,1] $, we obtain
  \begin{equation}\tag{A.3}
  \label{ReprFunc}
    T^{-3/2} \matW_{\trunc{Ts}} \matX_{\trunc{Ts}}' \bfeps_{\trunc{Ts}}
    = \tau( T^{-1/2} \epsilon_{\trunc{T \cdot }} )(s) + o_P(1),
  \end{equation}
  the $ o_P(1) $ being uniform in $ s \in [\kappa,1] $. It is easy
  to see that for any sequence $ \{ z, z_n \} \subset D_{\R^p}[\kappa,1] $ such that
  $ \lim_{n \to \infty} d_p(z_n,z) = 0 $, as $ n \to \infty $, with $ z \in C[\kappa,1] $, we have
  \[
    \lim_{n \to \infty} d_p(\tau( z_n ), \tau(z) ) = 0.
  \]
  Thus, the continuous mapping theorem in general separable metric
  spaces (Shorack and Wellner (1986), Th. 4, p. 47, and Remark 2, p.
  49) and (E) yield the result. $ \hfill \Box $
\end{myproofa}

\noindent
\begin{myproofa}{of Lemma~\ref{Lemma2}}{2}
Let $ \mbox{cond}_2(\matA(x)) = \sigma_k(x) / \sigma_1(x) $ denote
the condition of $ \matA(x) $ w.r.t. the spectral vector norm $ \|
\cdot \|_2 $. Let $ \varepsilon > 0 $. If $ \| \matA_n(x) -
\matA(x) \|_2 < \varepsilon $, the a-priori error estimate for
linear equations with disturbed coefficient matrices yields
\[
  \| a_{nj}^{-1}(x) - a_j^{-1} \|_2 \le
  \frac{ \mbox{cond}_2( \matA(x) ) }{ \| \matA(x) \|_2 - \varepsilon
  \mbox{cond}_2( \matA(x) ) } \varepsilon \| \veca_j(x)^{-1} \|,
\]
where $ a_{nj}^{-1}(x) $ ($ a_j^{-1}(x) $) denotes the $j$th
column of $ \matA_n^{-1}(x) $ ($ \matA^{-1}(x) $). $ \hfill \Box $
\end{myproofa}

\noindent
\begin{myproofa}{of Theorem~\ref{SeqRes1}}{3}
  Recall the representations
  \[
    \bhbeta_{\trunc{Ts}} - \bfbeta
      = ( \matX_{\trunc{Ts}}' \matX_{\trunc{Ts}} )^{-1} \matX_{\trunc{Ts}}' \bfeps_{\trunc{Ts}}
    \qquad \mbox{and} \qquad
    \widehat{\epsilon}_{\trunc{Ts}} = \epsilon_{\trunc{Ts}} - \vecx_{\trunc{Ts}}'( \bhbeta_{\trunc{Ts}} - \bfbeta
    ),
  \]
  where $ \bhbeta_{\trunc{Ts}} $ is defined in (\ref{DefBetaHat})
  yielding
  \begin{align*}
    \widehat{\epsilon}_{\trunc{Ts}}
      &= \epsilon_{\trunc{Ts}}
        - \vecx_{\trunc{Ts}}' (\matX_{\trunc{Ts}}'\matX_{\trunc{Ts}})^{-1} \matX_{\trunc{Ts}}' \bfeps_{\trunc{Ts}} \\
      &=\epsilon_{\trunc{Ts}}  -
        T^{1/2} \vecx_{\trunc{Ts}}' \matW_{\trunc{Ts}}
        ( \matW_{\trunc{Ts}} T^{-1}\matX_{\trunc{Ts}}'\matX_{\trunc{Ts}} \matW_{\trunc{Ts}} )^{-1}
        T^{-3/2} \matW_{\trunc{Ts}} \matX_{\trunc{Ts}}' \bfeps_{\trunc{Ts}}.
  \end{align*}
  Since
  $ \vecx_{\trunc{Ts}}' \matW_{\trunc{Ts}} = (1, \trunc{Ts}, \dots, \trunc{Ts}^p) \matW_{\trunc{Ts}} = \eins'$
  where $ \eins = (1,\dots,1)' \in \R^{p+1} $, we have
  \[
    T^{-1/2} \widehat{\epsilon}_{\trunc{Ts}}
    =
    T^{-1/2} \epsilon_{\trunc{Ts}} - \eins' \widetilde{\matH}_{\trunc{Ts}}^{-1}
    T^{-3/2} \matW_{\trunc{Ts}} \matX_{\trunc{Ts}}' \bfeps_{\trunc{Ts}}
  \]
  where
  \begin{equation}\tag{A.4}
  \label{DefHTs}
    \widetilde{\matH}_{\trunc{Ts}}
      = \matW_{\trunc{Ts}} T^{-1} \matX_{\trunc{Ts}}' \matX_{\trunc{Ts}} \matW_{\trunc{Ts}}.
  \end{equation}
  We will show that the matrix $ \widetilde{\matH}_{\trunc{Ts}} $,
  which equals
  \[
    \widetilde{\matH}_{\trunc{Ts}} = \left( \frac{\trunc{Ts}}{T} \trunc{Ts}^{-(i+j-1)} \sum_{t=1}^{\trunc{Ts}} t^{i+j-2} \right)_{i,j \in \{1,\dots,p\}},
  \]
  converges to $ s \matH $, uniformly in $ s \in [\kappa,1] $. For $ i = j = 1 $ this is
  obvious. Otherwise, $ i + j \ge 3 $, and for the corresponding elements of $ \widetilde{\matH}_{\trunc{Ts}} $ the assertion
  follows from
  \[
    \sum_{t=1}^{\trunc{Ts}} t^{i+j-2} = \frac{ (\trunc{Ts}+1)^{i+j-1} }{ i+j-1 }
    + O( (\trunc{Ts}+1)^{i+j-3} ).
  \]
  Hence
  \[
    \trunc{Ts}^{i+j-1} \sum_{t=1}^{\trunc{Ts}} t^{i+j-2}
      = \frac{1}{i+j-1} + O( 1 / (\trunc{T \kappa} + 1) ),
  \]
  yielding
  \begin{equation}\tag{A.5}
  \label{HilbertIsOP1}
    \sup_{s \in [\kappa,1]} | ( \widetilde{\matH}_{\trunc{Ts}} )_{ij} - s/(i+j-1) |
    = O( \trunc{T\kappa}^{-1} )
  \end{equation}
  for $ i,j \in \{ 1, \dots, p+1 \} $.
  Recall the representation (\ref{RepMatr}) and (\ref{ReprFunc}). Since $ \sup_{s \in
  [\kappa,1]} \tau( T^{-1/2} \bfeps_{\trunc{T \cdot}} ) $ converges
  weakly to the random variable
  \[
    \sup_{s \in [\kappa,1]} \tau( \eta B )(s) = \sup_{s \in
    [\kappa,1]} \biggl( \eta \int_0^s (r/s)^{i-1} B(r) \, dr
    \biggr)_{i=1,\dots, p+1},
  \]
  we may conclude that
  \begin{equation}\tag{A.6}
  \label{SupOP1}
    \sup_{s \in [\kappa,1]} T^{-3/2} \matW_{\trunc{Ts}}' \matX_{\trunc{Ts}}' \bfeps_{\trunc{Ts}} = O_P(1).
  \end{equation}
  (\ref{HilbertIsOP1}), (\ref{SupOP1}), and Lemma~\ref{Lemma2} imply
  \[
    \widehat{E}_T(s) =
      T^{-1/2} \epsilon_{\trunc{Ts}}
  - \eins' s^{-1} \matH^{-1} T^{-3/2} \matW_{\trunc{Ts}} \matX_{\trunc{Ts}}' \bfeps_{\trunc{Ts}}
  + O_P( 1 / \trunc{Ts} )
  \]
  Using the result (\ref{ReprFunc})
  we obtain
  \[
    \widehat{E}_T(s) =
      T^{-1/2} \epsilon_{\trunc{Ts}}
  - \eins' s^{-1} H^{-1} \tau( T^{-1/2} \bfeps_{\trunc{T \cdot}} ) + o_P(1)
  \]
  which shows that up to terms of order $ o_P(1) $ the process $ \widehat{E}_T $ is a continuous functional of
  $ \{ T^{-1/2} \epsilon_{\trunc{Ts}} : s \in [\kappa,1] \} $.
  Consequently,
  \[
    \widehat{E}_T(s)
    \stackrel{w}{\to} \eta \left\{ B(s) - s^{-1} \eins' H^{-1} \int_0^s (1,r/s,\dots,(r/s)^p)' B(r) \, dr
      \right\},
  \]
  in $ D[\kappa,1] $, as $ T \to \infty $. $ \hfill \Box $
\end{myproofa}

\noindent
\begin{myproofa}{of Theorem~\ref{FCLTUpdatedResiduals}}{4}
  The proof is similar as the proof of Theorem~\ref{SeqRes1}.
  We have
  \[
    \widehat{\epsilon}_{\trunc{Tr}}( \trunc{Ts} )
    =
    \epsilon_{\trunc{Tr}} - \vecx_{\trunc{Tr}}' \matW_{\trunc{Ts}}
    ( \matW_{\trunc{Ts}} \matX_{\trunc{Ts}}' \matX_{\trunc{Ts}} \matW_{\trunc{Ts}} )^{-1}
    \matW_{\trunc{Ts}} \matX_{\trunc{Ts}}' \bfeps_{\trunc{Ts}}.
  \]
  Note that
  \[
     \vecx_{\trunc{Tr}}' \matW_{\trunc{Ts}}
     = (1,\trunc{Tr}/\trunc{Ts},\dots,(\trunc{Tr}/\trunc{Ts})^p)'
  \]
  and let $ \vecv(r,s) = (1,r/s,\dots,(r/s)^p)' $. Due to (\ref{PowerFrac})
  we have
  \[
    \sup_{0 \le r \le s \le 1} \| \vecx_{\trunc{Tr}}'
    \matW_{\trunc{Ts}} - \vecv(r,s) \| = O(1/T).
  \]
  Combining this fact with (\ref{ReprFunc}) yields
  \begin{eqnarray*}
    T^{-1/2} \widehat{\epsilon}_{\trunc{Tr}}( \trunc{Ts} )
    &=&
    \frac{ \epsilon_{\trunc{Tr}} }{ T^{1/2} }
    -
    \vecx_{\trunc{Tr}}' \matW_{\trunc{Ts}} \widetilde{\matH}_{\trunc{Ts}}^{-1} \bfeps_{\trunc{Ts}} \\
    &=&
    \frac{ \epsilon_{\trunc{Tr}}( \trunc{Ts} ) }{ T^{1/2} }
    - \left\{ \vecv(r,s) s^{-1} \matH^{-1} + o_P(1) \right\}
      \left\{ \tau ( T^{-1/2} \epsilon_{\trunc{T \circ}} )(s) + o_P(1) \right\},
  \end{eqnarray*}
  where the $ o_P(1) $ terms are uniform in $ r,s \in [\kappa,1] $. Hence,
  uniformly in $ r,s \in [\kappa,1] $,
  $$
    T^{-1/2} \widehat{\epsilon}_{\trunc{Tr}}( \trunc{Ts} )
    = \varphi( T^{-1/2} \epsilon_{\trunc{T \cdot}} )(r,s) + o_P(1),
  $$
  where the functional $ \varphi : D[\kappa,1] \to D([\kappa,1]^2;R) $ is given by
  \begin{equation}\tag{A.7}
  \label{DefVarPhi}
    \varphi(z)(r,s) = z(r) - s^{-1} \vecv(r,s)  \matH^{-1} \int_0^s ( 1, u/s, \dots, (u/s)^p )'
    z(u) \, du, \qquad r,s \in [\kappa,1],
  \end{equation}
  for $ z \in D[\kappa,1] $. It is easy to see that for any
  sequence $ \{ z, z_n \} \subset D[\kappa,1] $ with $ d(z_n,z) \to 0 $, as $ n \to \infty $,
  and $ z \in C[\kappa,1] $, we have $ \| \varphi(z_n) - \varphi(z) \|_\infty \to 0 $,
  as $ n \to \infty $. Hence, an application of the continuous mapping theorem
  yields
  $$
    T^{-1/2} \widehat{\epsilon}_{\trunc{Tr}}( \trunc{Ts} )
    \stackrel{w}{\to}
     \varphi( \sigma B ) =
   \eta \left\{
    B(r) - s^{-1} \vecv(r,s)  \matH^{-1} \int_0^s (1,u/s,\dots, (u/s)^p)' B(u) \, du
  \right\},
  $$
  as $ T \to \infty $. $ \hfill \Box $
\end{myproofa}

\noindent
\begin{myproofa}{of Theorem~\ref{KPSS}}{5} We formulate the proof such
that the corresponding result for the change-point model of
Section~4 can be obtained by straightforward modifications. To
simplify exposition we assume $ \gamma = 0$. Note that for any $
\lambda_1, \lambda_2 \in \R $ the process
\[
  W_{\lambda_1,\lambda_2}(s) = \frac{\lambda_1}{\trunc{Ts}^4}
  \sum_{i=1}^{\trunc{Ts}} \bigl( \sum_{j=1}^i \widehat{\epsilon}_j(\trunc{Ts}) \bigr)^2 K(
  (i-\trunc{Ts})/h ) + \frac{\lambda_2}{\trunc{Ts}^2}
  \sum_{i=1}^{\trunc{Ts}} \widehat{\epsilon}_i^2(\trunc{Ts})
\]
can be written as
\begin{align*}
\tau_{\lambda_1,\lambda_2}( \widehat{E} ) & = \lambda_1 \biggl(
\frac{T}{\trunc{Ts}} \biggr)^4 \int_0^s \biggl( \int_0^r
\widehat{E}_{\trunc{Tz}}( \trunc{Ts} ) \, dz \biggr)^2 K( (\trunc{Tr}-\trunc{Ts})/h ) \, dr  \\
 & \qquad +
\lambda_2 \biggl( \frac{T}{\trunc{Ts}} \biggr)^2 \int_0^s
\widehat{E}_{\trunc{Tr}}^2(\trunc{Ts}) \, dr,
\end{align*}
where $ \tau_{\lambda_1,\lambda_2} $ maps elements of $ D([0,1]^2)
$ to elements of $ D[0,1] $. Let us show continuity of
$\tau_{\lambda_1,\lambda_2} $ w.r.t. the supnorm. W.l.o.g. we may
assume $ \| K \|_\infty = 1 $. Using the inequality $ |a^2 - b^2|
\le (|a|+|b|)|a-b| $ for real numbers $ a,b $ we can bound $
  | \tau_{\lambda_1,\lambda_2}(z_1) - \tau_{\lambda_1,\lambda_2}(z_2) |
$
by
\[
  [ \lambda_1 (T/\trunc{Ts})^4 \| z_1 \|_\infty + \lambda_2
  (T/\trunc{Ts})^2 (\| z_1 \|_\infty + \| z_2 \|_\infty) ] \| z_1 - z_2
  \|_\infty = O( \| z_1 - z_2 \|_\infty ).
\]
Hence, for $ 0 \le s_1 \le \cdots \le s_L \le 1 $, $ L \in \N $,
any associated linear combination $ \sum_{k=1}^L \rho_k
W_{\lambda_1,\lambda_2}(s_k) $, $ \rho_1, \dots, \rho_L \in \R $,
of the coordinates of the random vectors $
(W_{\lambda_1,\lambda_2}(s_1), \dots, W_{\lambda_1,\lambda_2}(s_L)
) $, converges in distribution to $ \sum_{k=1}^L \rho_k
\tau_{\lambda_1,\lambda_2}( \calE )(s_k) $, since $ \widehat{E}_T
\stackrel{w}{\to} \calE $, $ T \to \infty $, by
Theorem~\ref{FCLTUpdatedResiduals}. This verifies convergence of
the finite-dimensional distributions of the
$(D[\kappa,1])^2$-valued stochastic process $ (Z_{T1}, Z_{T2}) $,
where
\[
  Z_{T1}(s) =
    \trunc{Ts}^{-4} \sum_{i=1}^{\trunc{Ts}} \biggl(
    \sum_{j=1}^i \widehat{\epsilon}_j(\trunc{Ts}) \biggr)^2
    K((i-\trunc{Ts})/h), \quad
  Z_{T2}(s) = \trunc{Ts}^{-2} \sum_{i=1}^{\trunc{Ts}}
    \widehat{\epsilon}_i(\trunc{Ts})
\]
for $ s \in [\kappa,1] $. Tightness w.r.t. the product topology is
a consequence of Prohorov's theorem, since both coordinate
processes converge weakly. Thus, $ (Z_{T1}, Z_{T2})
\stackrel{w}{\to} (Z_1, Z_2) $, $ T \to \infty $, in $
(D[\kappa,1])^2 $, where
\[
  Z_1(s)= s^{-4} \int_0^s \bigl( \int_0^r \calE(z,s)
  \, dz \bigr)^2 K( \zeta(s-r) ) \, dr , \quad
  Z_2(s) = s^{-2} \int_0^s \calE^2(r,s) \,
  dr,
\]
for $ s \in [\kappa,1] $. Now a straightforward argument implies
that the ratio, $ V_T $, converges weakly, as $ T \to \infty $.
Finally, by Lipschitz continuity of the kernel $K$ the limit
process $ \calV $ is continuous w.p. $1$. $ \hfill \Box $
\end{myproofa}

\noindent
\begin{myproofa}{of Corollary~\ref{TheCorollary}}{6}
The result is shown using similar arguments as in Steland (2007b,
Corollary~2.1.) $ \hfill \Box $
\end{myproofa}

\section*{APPENDIX B: PROOFS OF RESULTS FROM SECTION~4}

By virtue of the method of proof used in the previous section, we
are in a position to extend the results for the kernel weighted
variance ratio process and its associated stopping time to the
change-point model of Section~3, if we have a FCLT for the process
of sequentially updated residuals. Thus we provide a detailed
proof of Theorem~\ref{FCLTUpdatedResidualsCP} and indicate the
required modifications to prove Theorem~\ref{FCLTCP}.

\noindent
\begin{myproofb}{of Theorem~\ref{FCLTUpdatedResidualsCP}}{1}
 Under the change-point model we have
 \[
 \{ T^{-1/2} \epsilon_{\trunc{Ts}} : \kappa \le s < \vartheta \}
 \stackrel{w}{\to}
  \{ \eta B(s) : \kappa \le s < \vartheta \},
 \]
 as $ T \to \infty $. Consider the process $ T^{-1/2} \epsilon_{\trunc{Ts}} $ for  $ \vartheta \le s \le 1
 $. First note that
 \[
   T^{-1/2} \epsilon_{\trunc{Ts}} \le \sup_{z \in [\kappa,1]} |
   T^{1/2} \epsilon_{\trunc{Tz}} |.
 \]
 Let $ \delta > 0 $.  By assumptions (C1) and (C2)
 \begin{align*}
  P \biggl( \sup_{z \in [\vartheta,1]} | T^{-1/2}
  \epsilon_{\trunc{Tz}} | > \delta \biggr)
  & = P \left( \max_{t = \trunc{T\vartheta}, \dots, T} | u_t | >
  T^{1/2} \delta / \xi_T \right) \\
  & \le (T - \trunc{T \vartheta} + 1)
          P( | u_1 | > T^{1/2} \delta / \xi_T ) \\
  & = O( T^{1-\gamma(1/2-\beta)} ) \\
  & = o_P( 1 ).
 \end{align*}
 if $ \beta < 1/2 $ and $ \gamma > (1/2-\beta)^{-1} $. Again using
 the Skorohod-Dudley-Wichura representation theorem we may assume that
 \[
   \sup_{s \in [\kappa,\vartheta]} | T^{-1/2}
   \epsilon_{\trunc{Ts}} - \eta B(s) | \stackrel{a.s.}{\to} 0,
 \]
 and
 \[
   \sup_{s \in [\kappa, \vartheta]} | T^{-1/2}
   \epsilon_{\trunc{Ts}} | \stackrel{a.s.}{\to} 0,
 \]
 as $ T \to \infty $. This implies a.s. convergence in the
 Skorohod metric to the cadlag process $ B
 \eins_{[\kappa,\vartheta)} $, i.e.,
 \[
  d( T^{-1/2} \epsilon_{\trunc{T \cdot}}, \eta B
  \eins_{[\kappa,\vartheta)} ) \stackrel{a.s.}{\to} 0,
 \]
 as $ T \to \infty $, which in turn yields weak convergence,
 \[
   T^{-1/2} \epsilon_{\trunc{Ts}} \stackrel{w}{\to} \eta B(s)
   \eins_{[\kappa,\vartheta)}(s),
 \]
 in $ D[\kappa,1] $, as $ T \to \infty $. Combining this fact with
 (\ref{ReprFunc}), the continuity of the functional $ \tau $
 (Jacod and Shiryaev~(2003), VI, Proposition~1.22, p.~329) yields
 \[
   T^{-3/2} \matW_{\trunc{Ts}} \matX_{\trunc{Ts}}'
   \bfeps_{\trunc{Ts}}
   \stackrel{w}{\to} \eta \int_0^{s} (1, u/s, \dots,
   (u/s)^p )' B(r) \, dr \eins_{[\kappa,\vartheta)}(s),
 \]
 as $ T \to \infty $. The same arguments as in the proof of Theorem~\ref{FCLTUpdatedResiduals}
 show that
 \[
   \widehat{E}_{\trunc{Tr}}( \trunc{Ts} )
   = T^{-1/2} \widehat{\epsilon}_{\trunc{Tr}}( \trunc{Ts} )
   = \varphi( T^{-1/2} \epsilon_{\trunc{T \cdot }} )(r,s) + o_P(1),
 \]
 as $ T \to \infty $, where the functional $ \varphi $ is defined in
 (\ref{DefVarPhi}). We have by linearity
 \begin{align*}
   \varphi( T^{-1/2} \epsilon_{\trunc{T \cdot }} )(r,s)
   & = \biggl[ T^{-1/2} \epsilon_{\trunc{Ts}} - s^{-1} \vecv(r,s) \matH^{-1} \int_0^s \vecv(u,s) T^{-1/2} \epsilon_{\trunc{Tu}} \, du
   \biggr]\\
    & \, = \biggl[ \eta B(s) - s^{-1} \vecv(r,s) \matH^{-1} \int_0^s \vecv(u,s) \eta B(u) \, du \biggr]
   \eins_{[\kappa,\vartheta)}(s) \\
   & \qquad + R_{T1}(s) -  R_{T2}(r,s) \\
   & \, = \varphi( \eta B \eins_{[\kappa,\vartheta)} )(r,s) +
   R_1(s) + R_2(r,s),
 \end{align*}
 where the remainder terms are given by
 \begin{align*}
   R_{T1}(s)   &=  T^{-1/2} \epsilon_{\trunc{Ts}} - \eta B(s)
   \eins_{[\kappa,\vartheta)}(s),\\
   R_{T2}(r,s) &=  s^{-1} \vecv(r,s) \matH^{-1} \int_0^s \vecv(u,s)[ T^{-1/2} \epsilon_{\trunc{Tu}} - \eta B(u)
   ] \eins_{[\kappa,\vartheta)}(s) \, du.
 \end{align*}
 Clearly, $  \sup_{s \in [\kappa,\vartheta)} | R_{T1}(s) | \to 0 $, as $ T \to \infty $, a.s.
 To estimate $ R_{T2} $, denote the maximum vector norm and the induced matrix norm
 by $ \| \circ \|_\infty $ and observe that
 \begin{align*}
   & \biggl\| \int_0^s \vecv(u,s) [ T^{-1/2} \epsilon_{\trunc{Tu}} -
   \eta
   B(u) ] \eins_{[\kappa,\vartheta)}(s) \, du \biggr\|_\infty \\
   & \qquad \le \int_0^s \| \vecv \|_\infty \sup_{z \in
   [\kappa,\vartheta)} | T^{-1/2} \epsilon_{\trunc{Tz}} - \eta
   B(z) | \eins_{[\kappa,\vartheta)}(s) \, du \\
   & \qquad \le \| \vecv \|_\infty \vartheta \sup_{z \in
   [\kappa,\vartheta)} | T^{-1/2} \epsilon_{\trunc{Tz}} - \eta
   B(z) | \stackrel{a.s.}{\to} 0,
 \end{align*}
 where $ \| \vecv \|_\infty = \sup_{r,s \in [\kappa,\vartheta)} \| \vecv(r,s) \|_\infty < \infty
 $. Hence,
 \[
   \sup_{s \in [\kappa,\vartheta)} | R_{T2}(r,s) | \le
   \kappa^{-1} \| \vecv \|_\infty^2  \| \matH^{-1} \|_\infty \vartheta
   \sup_{z \in [\kappa,\vartheta]} | T^{-1/2} \epsilon_{\trunc{Tz}} - \eta B(z)| \, du
   \stackrel{a.s.}{\to} 0,
 \]
 as $ T \to \infty $. Consequently,
 \[
   \varphi( T^{-1/2} \epsilon_{\trunc{T \cdot }} )(r,s)
   \stackrel{a.s.}{\to} \varphi( \eta B \eins_{[\kappa,\vartheta)}
   ),
 \]
 as $ T \to \infty $, which implies via
 \[
   d( \varphi( \epsilon_{\trunc{T \cdot}} ), \varphi( \eta B
   \eins_{[\kappa,\vartheta)} ) ) \stackrel{a.s.}{\to} 0,
 \]
 as $ T \to \infty $, weak convergence which completes the proof. $ \hfill \Box $
\end{myproofb}

\noindent
\begin{myproofb}{of Theorem~\ref{KPSS_CP}}{2}
To proof goes along the lines of the proof of Theorem~\ref{KPSS}.
Notice that now the linear combinations $ \sum_{k=1}^L \rho_k
W_{\lambda_1, \lambda_2}( s_k ) $ converge weakly in distribution
to $ \sum_{k=1}^L \rho_k \tau_{\lambda_1,\lambda_2}(
\calE_\vartheta )(s_k) $. $ \hfill \Box $
\end{myproofb}

\end{document}